\documentclass[a4paper]{jpconf}
\usepackage{bm,color}
\usepackage{graphicx}
\usepackage{amsmath}

\begin{document}
\title{Operator algebra as an application of logarithmic representation of infinitesimal generators}

\author{Yoritaka Iwata}

\address{Institute of Innovative Research, Tokyo Institute of Technology, Tokyo, Japan\\
Department of Mathematics, Shibaura Institute of Technology, Tokyo, Japan}

\ead{iwata\_phys@08.alumni.u-tokyo.ac.jp}

\begin{abstract}
The operator algebra is introduced based on the framework of logarithmic representation of infinitesimal generators.
In conclusion  a set of generally-unbounded infinitesimal generators is characterized as a module over the Banach algebra.
\end{abstract}

\section{Introduction}
The logarithmic representation of infinitesimal generator is introduced in Ref.~\cite{Iw-1}.
Its role is further discussed in Ref.~\cite{Iw-2}.
Let $X$ be a Banach space, and $Y$ be a dense Banach subspace of $X$.
Let $\{ U(t,s) \}_{-T \le t,s \le T}$ be an evolution family of operators satisfying the semigroup property (for example, see \cite{79tanabe}):
\[ \begin{array}{ll}
U(t,r) U(r,s) = U(t,s),  \vspace{1.5mm} \\
U(s,s) = U(s,t) U(t,s) = I
\end{array}  \]
on $X$.
Let the infinitesimal generator of $U(t,s)$ be denoted by $A(t)$ from $Y$ to $X$.
The logarithmic representation for the infinitesimal generator is
\begin{equation} \label{logex} \begin{array}{ll}
 A(t) ~ u_s =　
   (I+ \kappa U(s,t))~ \partial_{t} {\rm Log} ~ (U(t,s) + \kappa I) ~ u_s,
\end{array} \end{equation}
where $u_s$ is an element in $Y$, $\kappa \ne 0$ is a certain complex number, and ${{\rm Log}}$ denotes a principal branch of the logarithm. 
$\partial_t$ stands for a kind of weak differential defined in Ref.~\cite{Iw-1}.
For this definition, operators $U(t,s)$ and $A(t)$ are assumed to commute~\cite{Iw-1}.
In particular the representation plays a role of extracting the bounded part of the infinitesimal generator.
In this article a formulation of operator algebra is established for generally-unbounded infinitesimal generators.

\section{Basic relation}
The relation between exponential and logarithm functions of operators is shown under the validity of logarithmic representation (\ref{logex}).
According to Eq.~(\ref{logex}) the logarithm function in this case is given by $ {\rm Log} ~ (U(t,s) + \kappa I) $, and is formally or expectantly equal to the indefinite integral 
\[ \begin{array}{ll}
 {\rm Log} ~  (U(t,s) + \kappa I) =  {\displaystyle \int  (I+ \kappa U(s,t))^{-1} A(t) ~ dt }.
\end{array} \]
Its exponential function can be calculated as 
\begin{equation} \label{boundexp} \begin{array}{ll}
 \exp [ {\rm Log} ~  (U(t,s) + \kappa I)] =
 　 U(t,s) + \kappa I
\end{array} \end{equation}
${\rm Log}   \exp [ {\rm Log} ~  (U(t,s) + \kappa I)] =　{\rm Log} (U(t,s)+\kappa I)$ is always satisfied, but the above formal equality $ \int  (I+ \kappa U(s,t))^{-1} A(t) ~ dt  =  {\rm Log} ~  (U(t,s) + \kappa I) $ is not necessarily satisfied because the limited range of imaginary spectral distribution is true only for the right hand side.
In this sense ${\rm Log} ~ (U(t,s) + \kappa I)$ corresponds to the extracted bounded part of the infinitesimal generator $A(t)$.

A set of these operators is a candidate for bounded-operator algebra. 
Indeed, if $\kappa$ is taken from the resolvent set of operator $U(t,s)$, $ U(t,s) + \kappa I$ and  $(U(t,s) + \kappa I)^{-1}$ are necessarily bounded on $X$, and then the boundedness of $ {\rm Log} ~  (U(t,s) + \kappa I)$ follows.
That is,  $ {\rm Log} ~  (U(t,s) + \kappa I)$ and $U(t,s) + \kappa I$ are related by the exponential function of operators which can be defined by the convergent power series.
Meanwhile the logarithm function ${\rm Log} ~  (U(t,s) + \kappa I)$ cannot necessarily be defined by the convergent power series which is too restrictive to represent the function of $U(t,s)$ (i.e., $|U(t,s)-1|<1$ is required).  

\section{Main results}
Let ${\rm Log}(U(t,s)+\kappa I)$ be denoted by
\begin{equation} \label{ats} \begin{array}{ll}
a(t,s) = {\rm Log} ~ (U(t,s) + \kappa I),
\end{array} \end{equation}
where some algebraic properties of $a(t,s)$ can be found in Ref.~\cite{Iw-2}.  
Since $a(t,s)$ is bounded on $X$, $e^{a(t,s)}$ is well-defined by the convergent power series.
Consequently the exponentiability is reduced to ${\rm Log} ~ (U(t,s) + \kappa I)$.  
Using Eq.~(\ref{boundexp}),
\[ \begin{array}{ll}
U(t,s) = e^{a(t,s)} - \kappa I,
\end{array} \]
follows for a certain complex number $\kappa$.
It is understood by this equation that, with respect to the mathematical property of infinitesimal generators represented by Eq.~(\ref{logex}), it is sufficient to consider $a(t,s)$ instead of $A(t)$. \vspace{2.5mm}\\

\hspace{-6mm}
{\bf Theorem 1.}
Let $U_i(t,s)$ be evolution operators satisfying Eq.~(\ref{logex}), and ${\rm Log} ~ U_i(t,s)$ be well-defined for any $t,s \in [-T,T]$ and $i = 1,2, \cdots, n$.
${\rm Log} U_i(t,s)$ are assumed to commute with each other.
\[ \begin{array}{ll}
V_{Lg}(X) := \{ k {\rm Log} ~ U_i(t,s); ~ k \in {\bm C}, ~  t,s \in [-T,T] \} ~ \subset B(X)
    \end{array} \]
is a normed vector space over the complex number field, where $B(X)$ denotes a set of all the bounded operators on $X$, and the operator norm is equipped with $B(X)$. \vspace{2.5mm}\\

\hspace{-6mm}
{\bf Theorem 2.}
Let $U_i(t,s)$ be evolution operators satisfying Eq.~(\ref{logex}) for any $t,s \in [-T,T]$ and $i = 1,2, \cdots, n$.
For a certain $K \in B(X)$, let a subset of $B(X)$ in which each element is assumed to commute with ${\rm Log} ~ (U_i(t,s)+K)$ be $B_{ab}(X)$.
${\rm Log} (U_i(t,s)+K)$ are assumed to commute with each other.
\[ \begin{array}{ll}
B_{Lg}(X) := \left\{ {\mathcal K} {\rm Log} ~ (U_i(t,s) + K); ~  {\mathcal K} \in B_{ab}(X),~ K \in B(X) ,~  t,s \in [-T,T] \right\} ~ \subset B(X)
\end{array} \]
is a module over the Banach algebra $B(X)$. \vspace{2.5mm}\\

\section{Outlines of the proofs}
\subsection{A normed vector space}
The proof of Theorem 1 is presented.
In case of $\kappa = 0$ the operator $a(t,s)$ is reduced to
\[ \begin{array}{ll}
 {\rm Log} ~ U(t,s) \in B(X).
\end{array} \]
The operator sum is calculated using Dunford-Riesz integral
\begin{equation} \label{sum1} \begin{array}{ll}
  {\rm Log} ~ U(t,r) + {\rm Log} ~ U(r,s)   \vspace{2.5mm} \\
   = \frac{1}{2 \pi i} \int_{\Gamma}  {\rm Log} \lambda ~(\lambda - U(t,r))^{-1} d \lambda ~
    + \frac{1}{2 \pi i} \int_{\Gamma'}  {\rm Log} \lambda' ~(\lambda' - U(r,s))^{-1} d \lambda'  \vspace{2.5mm} \\
   = \frac{1}{(2 \pi i)^2}  \int_{\Gamma}  \int_{\Gamma'} ( {\rm Log} \lambda + {\rm Log} \lambda') 
     ~(\lambda - U(t,r))^{-1} ~(\lambda' - U(r,s))^{-1}
     ~ d \lambda'  d \lambda  \vspace{2.5mm} \\
   = \frac{1}{(2 \pi i)^2}  \int_{\Gamma}  \int_{\Gamma'} ( {\rm Log} \lambda  \lambda') 
     ~(\lambda - U(t,r))^{-1} ~(\lambda' - U(r,s))^{-1}
     ~ d \lambda'  d \lambda   \vspace{2.5mm} \\
   =   {\rm Log} ~ [ U(t,r) U(r,s)]
   =   {\rm Log} ~  U(t,s),
    \end{array} \end{equation}
then the sum closedness is clear.
Here $\Gamma'$ is assumed to be included in $\Gamma$, and this condition is not so restrictive in the present setting.
In a different situation, when $U(t,r)$ and $V(t,r)$ commute for the same $t$ and $r$, another kind of sum is calculated as
\begin{equation} \label{sum2} \begin{array}{ll}
  {\rm Log} ~ U_1(t,r) + {\rm Log} ~ U_2(t,r)   \vspace{2.5mm} \\
   = \frac{1}{2 \pi i} \int_{\Gamma}  {\rm Log} \lambda ~(\lambda - U_1(t,r))^{-1} d \lambda ~
    + \frac{1}{2 \pi i} \int_{\Gamma'}  {\rm Log} \lambda' ~(\lambda' - U_2(t,r))^{-1} d \lambda'  \vspace{2.5mm} \\
   = \frac{1}{(2 \pi i)^2}  \int_{\Gamma}  \int_{\Gamma'} ( {\rm Log} \lambda + {\rm Log} \lambda') 
     ~(\lambda - U_1(t,r))^{-1} ~(\lambda' - U_2(t,r))^{-1}
     ~ d \lambda'  d \lambda  \vspace{2.5mm} \\
   = \frac{1}{(2 \pi i)^2}  \int_{\Gamma}  \int_{\Gamma'} ( {\rm Log} \lambda  \lambda') 
     ~(\lambda - U_1(t,r))^{-1} ~(\lambda' - U_2(t,r))^{-1}
     ~ d \lambda'  d \lambda   \vspace{2.5mm} \\
   =   {\rm Log} ~ [ U_1(t,r) U_2(t,r)],
    \end{array} \end{equation}
where, for $W(t,r) = U_1(t,r) U_2(t,r)$, the semigroup property is satisfied as
\[ \begin{array}{ll}
 W(t,r) W(r,s) = U_1(t,r) U_2(t,r) U_1(r,s) U_2(r,s)  = U_1(t,r) U_1(r,s) U_2(t,r)  U_2(r,s) = W(t,s), \vspace{2.5mm}\\
 W(s,s) = W(s,t) W(t,s) = I,
    \end{array} \]
and then the sum closedness is clear.
Although the logarithm function is inherently a multi-valued function, the uniqueness of sum operation is ensured by the single-valued property of the principal branch ``Log''.
Consequently, since the closedness for scaler product is obvious,
\[ \begin{array}{ll}
V_{Lg}(X) = \{ k {\rm Log} ~ U(t,s); ~ k \in {\bf C}, ~  t,s \in [-T,T] \} \subset B(X)
    \end{array} \]
is a normed vector space over the complex number field.
In particular the zero operator ${\rm Log} I$ is included in $V_{Lg}(X)$.
Theorem 1 has been proved.

\subsection{$B(X)$-module} 
The proof of Theorem 2 is presented.
It is worth generalizing the above normed vector space.
In this sense, utilizing a common operator $K \in B(X)$, components are changed to $ {\rm Log} ~ (U_i(t,r)+ K)$.

The operator sum is calculated as
\begin{equation}  \begin{array}{ll}
  {\rm Log} ~ (U(t,r) + K) + {\rm Log} ~ (U(r,s)+ K)   \vspace{2.5mm} \\
   = \frac{1}{2 \pi i} \int_{\Gamma}  {\rm Log} \lambda ~(\lambda - U(t,r)-K )^{-1} d \lambda ~
    + \frac{1}{2 \pi i} \int_{\Gamma'}  {\rm Log} \lambda' ~(\lambda' - U(r,s)-K )^{-1} d \lambda'  \vspace{2.5mm} \\
   = \frac{1}{(2 \pi i)^2}  \int_{\Gamma}  \int_{\Gamma'} ( {\rm Log} \lambda + {\rm Log} \lambda') 
     ~(\lambda - U(t,r)-K)^{-1} ~(\lambda' - U(r,s)-K )^{-1}
     ~ d \lambda'  d \lambda  \vspace{2.5mm} \\
   = \frac{1}{(2 \pi i)^2}  \int_{\Gamma}  \int_{\Gamma'} ( {\rm Log} \lambda  \lambda') 
     ~(\lambda - U(t,r)- K)^{-1} ~(\lambda' - U(r,s)- K )^{-1}
     ~ d \lambda'  d \lambda   \vspace{2.5mm} \\
   =   {\rm Log} ~ [ (U(t,r) + K) (U(r,s) + K)]    \vspace{2.5mm} \\
   =   {\rm Log} ~ [ U(t,s) + K U(t,r) + K U(r,s) + K^2 ].
    \end{array} \end{equation}
After introducing a certain $K$ with sufficient large $|K|$, it is always possible to take integral path $\Gamma'$ to be included in $\Gamma$.
Since the part ``$K U(t,r) + K U(r,s) + K^2 I$'' is included in $B(X)$, the sum-closedness is clear. 
In a different situation, when $U_1(t,r)$ and $U_2(t,r)$ commute for the same $t$ and $r$, another kind of sum is calculated as
\begin{equation}  \begin{array}{ll}
  {\rm Log} ~ (U_1(t,r)+K) + {\rm Log} ~ (U_2(t,r)+K)   \vspace{2.5mm} \\
   = \frac{1}{2 \pi i} \int_{\Gamma}  {\rm Log} \lambda ~(\lambda - U_1(t,r)- K )^{-1} d \lambda ~
    + \frac{1}{2 \pi i} \int_{\Gamma'}  {\rm Log} \lambda' ~(\lambda' - U_2(t,r)-K )^{-1} d \lambda'  \vspace{2.5mm} \\
   = \frac{1}{(2 \pi i)^2}  \int_{\Gamma}  \int_{\Gamma'} ( {\rm Log} \lambda + {\rm Log} \lambda') 
     ~(\lambda - U_1(t,r)-K )^{-1} ~(\lambda' - U_2(t,r)-K )^{-1}
     ~ d \lambda'  d \lambda  \vspace{2.5mm} \\
   = \frac{1}{(2 \pi i)^2}  \int_{\Gamma}  \int_{\Gamma'} ( {\rm Log} \lambda  \lambda') 
     ~(\lambda - U_1(t,r)-K )^{-1} ~(\lambda' - U_2(t,r)-K )^{-1}
     ~ d \lambda'  d \lambda   \vspace{2.5mm} \\
  =   {\rm Log} ~ [ U_1(t,r) + K) (U_2(t,r) + K)]    \vspace{2.5mm} \\
  =   {\rm Log} ~ [ W(t,r) + K U_1(t,r) + K U_2(t,r) + K^2 I].
    \end{array} \end{equation}
Since the part ``$ K U_1(t,r) + K U_2(t,r) + K^2 $'' is included in $B(X)$, the sum-closedness is clear.

The product ${\mathcal K} {\rm Log} ~ (U_i(t,s) + K_2) \in B(X)$ is justified by the operator product equipped with $B(X)$. 
Consequently, since the closedness for operator product within $B(X)$ is obvious,
\[ \begin{array}{ll}
B_{Lg}(X) = \left\{ {\mathcal K} {\rm Log} ~ (U(t,s)+ K);  ~  {\mathcal K} \in B_{ab}(X),~ K \in B(X) ,~  t,s \in [-T,T] \right\}
\end{array} \]
is a module over the Banach algebra ($B(X)$-module).
In particular relation $V_{Lg}(X) \subset B_{Lg}(X)$ is satisfied.
Theorem 2 has been proved.

For the structure of $B_{Lg}(X)$, a certain originally unbounded part can be classified to ${\rm Log} ~ (U(t,s)+ K) \in B(X)$, and the rest part to ${\mathcal K} \in B_{ab}(X)$. 
Here the terminology ``originally unbounded'' is used, because some unbounded operators are reduced to bounded operators under the validity of the logarithmic representation.

\section*{Acknowledgement}
The author is grateful to Prof. Emeritus Hiroki Tanabe for fruitful comments.
This work was supported by grant-in-aid research (C) 17K05440. 　

\end{document}